\documentclass[a4paper,12pt]{article}
\usepackage[english]{babel}
\usepackage[T2A]{fontenc}
\usepackage[cp1251]{inputenc}
\usepackage{amsthm}
\usepackage[tbtags]{amsmath}
\usepackage{amsfonts,amssymb}
\begin{document}


\begin{center}
\textbf{Ekaterina Kompantseva, Askar Tuganbaev}

\textbf{Absolute Ideals\\ of Almost Completely Decomposable Abelian Groups}
\end{center}

\textbf{Abstract.} We consider the class $\mathcal{A}_0$ of Abelian block-rigid $CRQ$-groups of ring type. A subgroup $A$ of an Abelian group $G$ is called an \textsf{absolute ideal} of the group $G$ if $A$ is an ideal in any ring on $G$. We describe principal absolute ideals of groups in $\mathcal{A}_0$. This allows to prove that any group in $\mathcal{A}_0$ is an $afi$-group, i.e., a group $G$ such that any absolute ideal of $G$ is a fully invariant subgroup.


\textbf{MSC2020 datebase:} 16B99, 20K30, 20K99 

\section{Introduction}

Any homomorphism $\mu\colon G\otimes G\rightarrow G$ is called a \textsf{multiplication} on an Abelian group $G$. 
An Abelian group $G$ with multiplication defined on $G$ is called a \textsf{ring on} $G$. When studying rings on an abelian group $G$, one of the most natural questions is the question on subgroups of $G$ that are ideals in any ring on $G$. Such subgroups are called \textsf{absolute ideals} of the group $G$. 

Absolute ideals of Abelian groups were studied in \cite{Fuc73}, \cite{Che09}, \cite{Fri64}, \cite{Fri71}, \cite{Fuc66}, \cite{McL75b}.
In \cite{Fri64}, general properties of absolute ideals are studied. 
To do this, it is defined the group $\text{M}(G)$ generated by all homomorphic images of the group $G$ in its endomorphism group $\text{End }G$ and it is proved that a subgroup $A$ of an Abelian group $G$ is an absolute ideal of $G$ if and only if $A$ is $\text{M}(G)$-admissible, i.e., $\text{M}(G)(A)\subseteq A$. 
The results of this paper were generalized to \cite{Fri71} for modules over a commutative unital ring.
In \cite{McL75}, it is proved that any absolute ideal of an Abelian $p$-group $G$ is of the form $G=H\cap \{g\in G\,|\,p^ng\in S\}$, where $H$ is a fully invariant subgroup of the group $G$, $S$ is a subgroup of the group $G^1$, and $n$ is an integer. In \cite{Che09}, Abelian groups, in which all subgroups are absolute ideals, are described. It is proved that the class of all such groups consists of torsion groups, each $p$-component of which is a cyclic or divisible group, torsion-free nil-groups, as well as cyclic groups of infinite order.

We remark that for an Abelian group $G$, the subgroup $\text{M}(G)$ of $G$ is an ideal in the endomorphism ring $\text{E}(G)$. Therefore, there is a close connection between the completely invariant subgroups of an Abelian group and its absolute ideals. It is easy to see that any fully invariant subgroup of an Abelian group is its absolute ideal \cite{Fuc73}. However, the converse is not true. In \cite{Fri64}, the problem of describing Abelian groups in which any absolute ideal is a completely invariant subgroup is formulated. Such groups are called \textsf{$afi$-groups}. In \cite{Fri64}, it is proved that an Abelian group $G$ is an $afi$-group if and only if for each element $g\in G$, the subgroups $A(g)=\{\eta\in\text{End}G\,|\,\eta(g)=0\}$ and $\text{M}(G)$ together with the identity endomorphism generate the entire group $\text{End}G$. However, in reality, this criterion does not bring us closer to solving the problem, but is only a variant of the definition of an $afi$-group. In \cite{McL75b}, it is proved that a fully transitive $p$-group $G$ is an $afi$-group if and only if its first Ulm subgroup $G^1$ is a cocyclic group. Since any decomposition of a torsion group $G$ into the direct sum of its $p$-components also is a decomposition of any ring on $G$ into a direct sum of ideals, it is easy to obtain a description of fully transitive torsion $afi$-groups. The same result was obtained in \cite{Pha11}, and the methods used in the proof are essentially different from the methods in \cite{McL75b}. Mixed $afi$-groups are studied in \cite{Pha12} and \cite{KomP19} and torsion-free $afi$-groups are studied in \cite{KomT22b}.

In the study of absolute ideals, principal absolute ideals introduced in \cite{Pha11} play an important role. A \textsf{principal absolute ideal} of the Abelian group $G$ generated by an element $g\in G$ is the least  absolute ideal $\langle g\rangle_{AI}$ of the group $G$ containing $g$. Since any absolute ideal $K$ is a sum of principal absolute ideals, $K=\sum\limits_{g\in K}\langle g\rangle_{AI}$, the study of the properties of arbitrary absolute ideals is often reduced to the case of principal absolute ideals.

This paper is devoted to questions related to absolute ideals of almost completely decomposable Abelian groups. The paper is a continuation of the papers \cite{KT23} and \cite{KT24}, where the properties of the group $\text{Mult}\,G$ of all multiplications of an almost completely decomposable group $G$ were studied.

The paper considers only additively written Abelian groups and the word <<group>> everywhere means <<abelian group>>.

A torsion-free group $G$ of finite rank is called an \textsf{almost completely decomposable} group (an \textsf{$ACD$-group}) if $G$ contains a completely decomposable subgroup of finite index. $ACD$-groups are studied in \cite{Bla04}, \cite{Bla06}, \cite{Bla08}, \cite{BlaM94}, \cite{DugO93}, \cite{Kom09}, \cite{KT23}, \cite{KT24}, \cite{Mad00} and other papers. Every $ACD$-group $G$ contains a completely decomposable fully invariant subgroup $\text{Reg}\,G$ of finite index which is called the \textsf{regulator} of the group $G$. The factor group $G/\text{Reg}\,G$ is called the \textsf{regulator quotient} and the index $n(G)$ of the subgroup $\text{Reg}\,G$ in the group $G$ is called the \textsf{regulator index} of the group $G$. $ACD$-groups with cyclic regulator quotient are called \textsf{$CRQ$-groups}.

Let $G$ be an almost completely decomposable group. Then the group $\text{Reg}\,G$ is uniquely, up to isomorphism, represented as a direct sum of torsion-free groups of rank $1$ \cite[Proposition 86.1]{Fuc73}. For every type $\tau$, we denote by $\text{Reg}_{\tau}\,G$ the sum of direct summands of rank 1 and of type $\tau$ in this direct decomposition of the group $\text{Reg}\,G$. This sum is called the \textsf{$\tau$-homogeneous component} of the group $G$. 

The set of types
$$
T(G)=T(\text{Reg}\,G)=\{\tau\,|\,\text{Reg}_{\tau}\,G\ne 0\}
$$
is called the \textsf{set of critical types} of the groups $G$ and $\text{Reg}\,G$. 
If $T(G)$ consists of pairwise incomparable types, then the groups $G$ and $\text{Reg}\,G$ are called \textsf{block-rigid} groups. In this case, if for any $\tau\in T(G)$, the group $\text{Reg}_{\tau}G$ is of rank $1$, then $G$ and $\text{Reg}\,G$ are called \textsf{rigid} groups. If all types in $T(G)$ are idempotent types, then $G$ is called a group \textsf{of ring type}.

We denote by $\mathcal{A}_0$ the class of all reduced block-rigid $CRQ$-groups of ring type. In Section 2, we describe principal absolute ideals of groups in $\mathcal{A}_0$ (Theorem 2.4). This allows us to prove in Section 3 that any group in $\mathcal{A}_0$ is an $afi$-group (Theorem 3.2).

A multiplication $\mu\colon G\otimes G\rightarrow G$ is often denoted by $\times$, i.e., $\mu(g_1\otimes g_2)=g_1\times g_2$ for all $g_1,g_2\in G$. 
The ring on the group $G$, defined by this multiplication, is denoted by  $(G,\times)$.
Let $G$ be a group, $g\in G$, and let $(G,\times)$ be a ring on $G$. The rank of the group $G$ is denoted by $r(G)$, the divisible hull of $G$ is denoted by $\tilde{G}$. The endomorphism group and the endomorphism ring of the group $G$ are denoted by $\text{End}\,G$ and $E(G)$, respectively. We denote by $(g)_\times$ the principal ideal of the ring $(G,\times)$ generated by the element $g$. 
If $G$ is a completely decomposable group and $\tau$ is a type, then $G_{\tau}$ is the $\tau$-homogeneous component of the group $G$ and $\tau(p)$ is the component of the type $\tau$ corresponding to the integer $p$. 
If $S\subseteq G$, then $\langle S\rangle$ and $\langle S\rangle_*$ are the subgroup and the pure subgroup of the group $G$ generated by the set $S$, respectively. As usual, $\mathbb{N}$ and $P$ are the set of all positive integers and the set of all prime integers, respectively, $\mathbb{Z}$ is the ring of integers and $\mathbb{Q}$ is the field of rational numbers. If $S$ is a finite subset in $\mathbb{Z}$, then $\text{gcd}(S)$ is the greatest common divisor of all integers in $S$ and $\text{lcm}(A)$ is the the least common multiple of all integers in $S$.

For any type $\tau$, we set
$$
P_{\infty}(\tau)=\{p\in P \,|\, \tau(p)=\infty\}, \quad 
P_0(\tau)=P\setminus P_{\infty}(\tau).
$$
If $P_1\subseteq P$, then a \textsf{$P_1$-number} is a non-zero integer such that any its prime factor (if it exists) is contained in $P_1$.
If $R$ is a unital ring, then $Re$ is the cyclic $R$-module generated by the element $e$. If $\tau$ is an idempotent type, then $R_{\tau}$ is a unital subring of the field $\mathbb{Q}$ such that the additive group of $R_{\tau}$ is of type $\tau$. For notation and definitions, unless otherwise specified, we refer to \cite{Fuc15} and \cite{KMT03}.

\section{Principal Absolute Ideals of $CRQ$-Groups}

In this section, we describe principal absolute ideals of groups in the class $\mathcal{A}_0$.

In \cite{DugO93}, positive integers $m_{\tau}=m_{\tau}(G)$ are defined for a group $G\in \mathcal{A}_0$, and it is proved that the set $\{m_{\tau}\,|\,\tau\in T(G)\}$ is an invariant system of almost isomorphism of the group $G$. In addition, $n(G)=\text{lcm}\{m_{\tau}\,|\,\tau\in T(G)\}$. In \cite[Theorem 3.5]{BlaM94}, it is proved that for any group $G\in \mathcal{A}_0$, there exists a direct decomposition 
$$
G=G_1\oplus C,\eqno (2.1)
$$
where $C$ is a block-rigid completely decomposable group and $G'$ is a rigid $CRQ$-group satisfying the following conditions:
$$
\tau\in T(G_1) \text{ if and only if } m_{\tau}(G)> 1,\eqno (2.1')
$$
$$
m_{\tau}(G_1)=m_{\tau}(G) \text{ for all } \tau\in T(G').\eqno (2.1'')
$$

Decomposition $(2.1)$, satisfying conditions $(2.1')$ and $(2.1'')$, is called a \textsf{principal decomposition} of the group $G$. In a principal decomposition of $G$, the group $G'$ does not contain completely decomposable direct summands; such groups are said to be \textsf{clipped}. We remark that a principal decomposition of a $CRQ$-group is not uniquely defined \cite{BlaM94}. Let $\text{Reg}\,G_1=B$. Then 
$T(G_1)=T(B)=\{\tau\in T(G)\,|\,m_{\tau}>1\}$ and $\tilde{G}_1=\tilde{B}$.

Let $d$ be an element of $G_1$ such that $G/\text{Reg}\,G=\langle d+\text{Reg}\,G\rangle$. Then there exists a system $E_0=\{e_0^{(\tau)}\in B_{\tau}\,|\, \tau\in T(B)\}$ such that
$$
B=\oplus_{\tau\in T(B)}R_{\tau}e_0^{(\tau)}, \eqno (2.2)
$$

the element $d$ of the group $\tilde{B}$ can be represented in the form
$$
d=\sum_{\tau\in T(B)}\dfrac{s_{\tau}}{m_{\tau}}e_0^{(\tau)}, \eqno (2.3)
$$
where $m_{\tau}=m_{\tau}(G)$, $s_{\tau}\in \mathbb{Z}$, and the following conditions hold:
$$
\text{gcd}(s_{\tau},m_{\tau})=1 \text{ for all } \tau\in T(B).
\eqno (2.3')
$$
$$
s_{\tau} \text{ and } m_{\tau} \text{ are }P _0(\tau)\text{-numbers for every } \tau\in T(B).\eqno(2.3'')
$$

A system $E_0=\{e_0^{(\tau)}\in B_{\tau} \,|\, \tau\in T(B)\}$, satisfying conditions $(2.2)$ and $(2.3)$, is called an \textsf{$rc$-basis} of the group $G$ defined by the element $d$. We remark that the pair $(d,E_0)$ uniquely determines the numbers $s_{\tau}$ ($\tau\in T(B)$).
Equality $(2.3)$ is called a\\ 
\textsf{standard representation of a block-rigid $CRQ$-group $G$ related to the pair} $(d,E_0)$. 

In its turn, the group $C$ is of the form $C=\oplus_{\tau\in T(C)}C_{\tau}$, where $C_{\tau}$ can be represented in the form $C_{\tau}=\oplus_{i\in I_{\tau}(C)}R_{\tau}e_i^{(\tau)}$ and $I_{\tau}(C)$ is a non-empty finite subset in $\mathbb{N}$ for $\tau\in T(C)$. 
We set $I_{\tau}(B)=\{0\}$ for $\tau\in T(B)$. We also set $I_{\tau}(B)=\varnothing$ for $\tau\notin T(B)$ and $I_{\tau}(C)=\varnothing$ for $\tau\notin T(C)$. For any $\tau\in T(G)$, we set $I_{\tau}=I_{\tau}(B)\cup I_{\tau}(C)$. A system $E=\{e_i^{(\tau)}\,|\,\tau\in T(G),\;i\in I_{\tau}\}$ is called an \textsf{$r$-basis} of the group $G$ if its subsystem $E_0=\{e_0^{(\tau)}\,|\,\tau\in T(B)\}$ is an $rc$-basis of the group $G$.

We set $\text{Reg}\,G=A$. Then $A_{\tau}=B_{\tau}\oplus C_{\tau}$ and $A=\oplus_{\tau\in T(G)}A_{\tau}$. According to \cite[Proposition 2.4.11]{Mad00}, such a decomposition of a completely decomposable group is unique if and only if $A$ is a block-rigid group. Let $\tilde G$, $\tilde A$ and $\tilde A_{\tau}$ be divisible hulls of the groups $G$, $A$ and $A_{\tau}$, respectively. Then we have
$$
\tilde G=\tilde A=\oplus_{\tau\in T(G)}\tilde A_{\tau}.
$$

For $\tau\in T(G)$, we denote by $\pi_{\tau}$ the natural projection of the group $\tilde G$ on $\tilde A_{\tau}$. 

In what follows, we fix a principal decomposition and an $r$-basis of the group $G$.

Let $G\in \mathcal{A}_0$, $\text{Reg}\,G=A$, $E=\{e_i^{(\tau)}\,|\,\tau\in T(G),\;i\in I_{\tau}\}$ be an $r$-basis of the group $G$, and let $G$ have standard representation $(2.3)$. Since $G=\langle d,A\rangle$, any element of $G$ is of the form
$$
g=\sum_{\tau\in T(B)}\dfrac{r_{\tau}}{m_{\tau}}e_{0}^{(\tau)}+
\sum_{\tau\in T(C)}\sum_{i\in I_{\tau}(C)}a_{i}^{(\tau)}e_{i}^{(\tau)}\in G,
$$
where $r_{\tau},a_{i}^{(\tau)}\in R_{\tau}$ for $\tau\in T(G)$, $i\in I_{\tau}(C)$. 

\textbf{Remark 2.1.} Any set of elements $\{u_{ij}^{(\tau)}\in A_{\tau}\,|\,\tau\in T(G),\;i,j\in I_{\tau}\}$ determines a multiplication $\times$ on $A$ such that $e_{i}^{(\tau)}\times e_{j}^{(\tau)}=u_{ij}^{(\tau)}$ for any $\tau\in T(G)$ and each $i,j\in I_{\tau}$; in addition, $e_{i}^{(\tau)}\times e_{j}^{(\sigma)}=0$ for $\tau\ne \sigma$ and any $i\in I_{\tau}$, $j\in I_{\sigma}$. According to \cite[Theorem 2.4]{KT24}, this multiplication can be extended to a multiplication on $G$ if and only if there exists $\alpha\in\mathbb{Z}$ such that
$$
u_{i0}^{(\tau)}=m_{\tau}v_{i0}^{(\tau)},\;
u_{0i}^{(\tau)}=m_{\tau}v_{0i}^{(\tau)},\;
v_{00}^{(\tau)}=\alpha s_{\tau}^{-1}e_0+m_{\tau}a_{00}^{(\tau)}\eqno(2.4)
$$
for some $v_{i0}^{(\tau)},v_{0i}^{(\tau)},a_{00}^{(\tau)}\in A_{\tau}$. Here $s_{\tau}^{-1}$ is an integer which is converse to $s_{\tau}$ modulo $m_{\tau}$; such an integer exists by $(2.3')$.

If a ring $(G,\times)$ with multiplication $\times$ satisfies condition $(2.5)$, then we say that the ring \textsf{corresponds to the integer $\alpha$ with respect to the pair} $(d,E_0)$. 
The multiplication, corresponding to the integer $\alpha$, can be described if we consider it as an element of the group $\text{Mult }G$ of all multiplications of the group $G$. 
Let $M=\text{Mult }G$. In \cite{KT24}, it is proved that the group $M$ also belongs to the class $\mathcal{A}_0$ and it can be represented in the form $M=\langle X,\text{Reg }M\rangle$, where $X=X(d,E_0)$ is an element of the group $M$ depending on the pair $(d,E_0)$. The multiplication $\times$ on the group $G$ corresponds to $\alpha\in \mathbb{Z}$ with respect to the pair $(d,E_0)$ if and only if $\times$ belongs to the class $\alpha X+\text{Reg }M$.

The \textsf{principal absolute ideal} of the group $G$ generated by an element $g$ is the least absolute ideal $\langle g\rangle_{AI}$ of the group $G$ containing $g$. In \cite{Fri71}, it is defined the subgroup $M(G)=\langle\text{Im}\,\Psi\,|\,\Psi\in\text{Hom}\,(G,\text{End}\,G)\rangle$ of the group $\text{End}\,G$ and it is proved that $M(G)$ is an ideal of the endomorphism ring $E(G)$. We set $M(G)(g)=\langle \varphi(g)\,|\,\varphi\in M(G)\rangle$ for $g\in G$.

To describe principal absolute ideals of groups, we need the following lemma in \cite{Pha12}.

\textbf{Lemma 2.2 \cite{Pha12}.} Let $G$ be a group and $g\in G$. Then

\textbf{1)} $M(G)(g)=\langle g\times x\,|\,\times\in\text{Mult}\,(G),\,x\in G\rangle$,

\textbf{2)} $\langle g\rangle_{AI}=\langle g\rangle+M(G)(g)$.~$\blacktriangleright$

\textbf{Remark 2.3.} If $\tau$ is a type and $a\in R_{\tau}\setminus \{0\}$, then $a$ can be uniquely written in the form $a=\overline{a}y$, where $\overline{a}$ is a positive $P_0(\tau)$-number and $y$ is a $P_{\infty}(\tau)$-fraction. We define the number $\overline{a}$ for every $a\in R_{\tau}$, by setting $\overline{0}=0$.

\textbf{Theorem 2.4.} Let $G$ be a block-rigid $CRQ$-group of ring type with $r$-basis $E=\{e_i^{(\tau)}\,|\,\tau\in T(G),\;i\in I_{\tau}\}$,  standard representation $(2.3)$ and $\text{Reg }G=A$. Let
$$
g=\sum_{\tau\in T(B)}\dfrac{r_{\tau}}{m_{\tau}}e_{0}^{(\tau)}+
\sum_{\tau\in T(C)}\sum_{i\in I_{\tau}(C)}a_{i}^{(\tau)}e_{i}^{(\tau)}\in G,
$$
where $r_{\tau},a_{i}^{(\tau)}\in R_{\tau}$ for $\tau\in T(G)$, $i\in I_{\tau}(C)$. By setting $r_{\tau}=0$ for $\tau\notin T(B)$ and $a_{i}^{(\tau)}=0$ for $\tau\notin T(C)$ (in this case, we have $I_{\tau}(C)=\varnothing$), we set
$$
\ell_{\tau}(g)=\text{gcd}\,(\overline{r_{\tau}},\{\overline{a_{i}^{(\tau)}}\,|\,i\in I_{\tau}(C)\})\quad\text{ for }\,\tau\in T(G).
$$

Then 
$$
\langle g\rangle_{AI}=\langle g\rangle+\oplus_{\tau\in T(B)}r_{\tau}A_{\tau} +\oplus_{\tau\in T(C)}\oplus_{i\in I_{\tau}(C)}a_{i}^{(\tau)}A_{\tau}=
$$
$$
=\langle g\rangle+\oplus_{\tau\in T(G)}\ell_{\tau}(g)A_{\tau}.
$$

\textbf{Remark.} We can prove that for any $g\in G$ and every $\tau\in T(G)$, the number $\ell_{\tau}(g)$ does not depend on the choice of an $r$-basis $E$ and a standard representation of the group $G$.

$\blacktriangleleft$ Let $\tau\in T(G)$. We set
$$
L_{\tau}=r_{\tau}A_{\tau}+\oplus_{i\in I_{\tau}(C)}a_{i}^{(\tau)}A_{\tau},\qquad
L=\oplus_{\tau\in T(G)}L_{\tau},
$$

We prove that $\langle g\rangle_{AI}=\langle g\rangle+L$. 

First, we prove that $M(G)(g)\subseteq\langle g\rangle+L$. Let
$(G,\times)$ be a ring. Then there exists a ring $(\tilde{G},\times)$ on the divisible hull $\tilde{G}=\tilde{A}$ of the group $G$ and $G$ is a subring in $(\tilde{G},\times)$, \cite[Chapter 18, Theorem 1.3]{Fuc15}. In the ring $(G,\times)$, we set $e_{i}^{(\tau)}\times e_{j}^{(\tau)}=u_{ij}^{(\tau)}\in A_{\tau}$ for all $\tau\in T(G)$ and $i,j\in I_{\tau}$. It follows from Remark 2.1 that there exists an integer $\alpha$ such that for any $\tau\in T(B)$ and each $i\in I_{\tau}$, we have
$$
u_{i0}^{(\tau)}=m_{\tau}v_{i0}^{(\tau)},\;
u_{0i}^{(\tau)}=m_{\tau}v_{0i}^{(\tau)},\;
v_{00}^{(\tau)}=\alpha s_{\tau}^{-1}e_0+m_{\tau}a_{00}^{(\tau)}
$$
for some $v_{i0}^{(\tau)},v_{0i}^{(\tau)},a_{00}^{(\tau)}\in A_{\tau}$. Here $s_{\tau}^{-1}$ is an integer such that the relation $s_{\tau}s_{\tau}^{-1}=1+m_{\tau}y_{\tau}$ holds for some $y_{\tau}\in\mathbb{Z}$.

Let $\tau\in T(G)$ and $k\in I_{\tau}$. If $\tau\notin T(B)$, then 
$$
g\times e_{k}^{(\tau)}=\pi_{\tau}(g)\times e_{k}^{(\tau)}=
\left(\sum_{i\in I_{\tau}}a_{i}^{(\tau)}e_{i}^{(\tau)}\right)\times
e_{k}^{(\tau)}=\sum_{i\in I_{\tau}(C)}a_{i}^{(\tau)}u_{ik}^{(\tau)}\in L_{\tau}.
$$
In addition, $\pi_{\tau}(g\times d)=0$.

Let $\tau\in T(B)$. In the ring $(\tilde{G},\times)$, we have
$$
g\times e_{k}^{(\tau)}=\pi_{\tau}(g)\times e_{k}^{(\tau)}=
\left(\dfrac{r_{\tau}}{m_{\tau}}e_{0}^{(\tau)}+\sum_{i\in I_{\tau}(C)}a_{i}^{(\tau)}e_{i}^{(\tau)}\right)\times e_{k}^{(\tau)}=
$$
$$
=\dfrac{r_{\tau}}{m_{\tau}}u_{0k}^{(\tau)}+\sum_{i\in I_{\tau}(C)}a_{i}^{(\tau)}u_{ik}^{(\tau)}=
r_{\tau}v_{0k}^{(\tau)}+\sum_{i\in I_{\tau}(C)}a_{i}^{(\tau)}u_{ik}^{(\tau)}\in L_{\tau}.
$$
In addition, 
$$
\pi_{\tau}(g\times d)=\pi_{\tau}(g)\times \pi_{\tau}(d)=
\left(\dfrac{r_{\tau}}{m_{\tau}}e_{0}^{(\tau)}+
\sum_{i\in I_{\tau}(C)}a_{i}^{(\tau)}e_{i}^{(\tau)}\right)\times 
\dfrac{s_{\tau}}{m_{\tau}}e_{0}^{(\tau)}=
$$
$$
=\dfrac{r_{\tau}s_{\tau}}{m_{\tau}^2}u_{00}^{(\tau)}+
\sum_{i\in I_{\tau}(C)}\dfrac{a_{i}^{(\tau)}s_{\tau}}{m_{\tau}}u_{i0}^{(\tau)}=
\dfrac{r_{\tau}s_{\tau}\alpha s_{\tau}^{-1}}{m_{\tau}}e_{0}^{(\tau)}+r_{\tau}s_{\tau}a_{00}^{(\tau)}+
\sum_{i\in I_{\tau}(C)}a_{i}^{(\tau)}s_{\tau}v_{i0}^{(\tau)}=
$$
$$
=\alpha\dfrac{r_{\tau}}{m_{\tau}}(1+m_{\tau}y_{\tau})e_{0}^{(\tau)}+r_{\tau}s_{\tau}a_{00}^{(\tau)}+
\sum_{i\in I_{\tau}(C)}a_{i}^{(\tau)}s_{\tau}v_{i0}^{(\tau)}=
$$
$$
=\alpha\dfrac{r_{\tau}}{m_{\tau}}e_{0}^{(\tau)}+r_{\tau}(\alpha
y_{\tau}e_{0}^{(\tau)}+s_{\tau}a_{00}^{(\tau)})+
\sum_{i\in I_{\tau}(C)}a_{i}^{(\tau)}s_{\tau}v_{i0}^{(\tau)}=
$$
$$
=\alpha\left(\dfrac{r_{\tau}}{m_{\tau}}e_{0}+
\sum_{i\in I_{\tau}(C)}a_{i}^{(\tau)}e_{i}^{(\tau)}\right)-
\alpha\sum_{i\in I_{\tau}(C)}a_{i}^{(\tau)}e_{i}^{(\tau)}+
r_{\tau}b_{\tau} +\sum_{i\in I_{\tau}(C)}a_{i}^{(\tau)}s_{\tau} v_{i0}^{(\tau)}=
$$
$$
=\pi_{\tau}(\alpha g)+r_{\tau}b_{\tau}+
\sum_{i\in I_{\tau}(C)}a_{i}(s_{\tau}v_{i0}^{(\tau)}-\alpha e_{i}^{(\tau)})=
$$
$$
=\pi_{\tau}(\alpha g)+r_{\tau}b_{\tau}+
\sum_{i\in I_{\tau}(C)}a_{i}c_{i}^{(\tau)},
$$
where $b_{\tau}=\alpha y_{\tau}e_{0}^{(\tau)}+ s_{\tau}a_{00}^{(\tau)}\in A_{\tau}$, $c_{i}^{(\tau)}=s_{\tau}v_{i0}^{(\tau)}-\alpha e_{i}^{(\tau)}\in A_{\tau}$ for $i\in I_{\tau}(C)$. Therefore, $\pi_{\tau}(g\times d)\in \pi_{\tau}(\alpha g)+L_{\tau}$, whence $g\times d\in \langle g\rangle+L$. 

Since any element $x\in G$ can be represented in the form
$$
x=\beta d+\sum_{\tau\in T(G)}\sum_{k\in I_{\tau}}x_{k}^{(\tau)}e_{k}^{(\tau)},
$$

where $\beta \in\mathbb{Z}$, $x_{k}^{(\tau)}\in R_{\tau}$,
we have that $g\times x\in \langle g\rangle+L$ for any $x\in G$. Since the multiplication $\times$ is arbitrary, $M(G)(g)\subseteq \langle g\rangle+L$, whence
$$
\langle g\rangle_{AI}=\langle g\rangle+M(G)(g)\subseteq \langle g\rangle+L,
$$
by Lemma 2.2.

Now we prove that $L\subseteq \langle g\rangle_{AI}$. Let $\tau\in T(G)$. We prove that $L_{\tau}\subseteq\langle g \rangle_{AI}$. Let $k\in I_{\tau}$, $b\in R_{\tau}$.

First, we assume that $\tau\in T(B)$ and prove that $r_{\tau}be_{k}^{(\tau)}\in \langle g \rangle_{AI}$.

\textbf{Case 1.} $\tau\in T(B)\cap T(C)$.

We fix $t\in I_{\tau}(C)$. According to Remark 2.1, relations
$$
e_{0}^{(\tau)}\times e_{t}^{(\tau)}=e_{t}^{(\tau)}\times e_{0}^{(\tau)}=m_{\tau}e_{k}^{(\tau)},
$$
$e_{i}^{(\tau)}\times e_{j}^{(\sigma)}=0$ if $\sigma\ne \tau$ or $(i,j)\ne (0,t)$, $(i,j)\ne (t,0)$, 

define a multiplication $\times$ on $G$. We consider the ring $(G,\times)$ as a subring of the ring $(\tilde{G},\times)$ \cite[Chapter 18, Theorem 1.3]{Fuc15} and obtain
$$
g\times be_{t}^{(\tau)}=\dfrac{r_{\tau}}{m_{\tau}}e_{0}^{(\tau)}\times be_{t}^{(\tau)}=r_{\tau}be_{k}^{(\tau)}\in \langle g \rangle_{\times}\subseteq\langle g \rangle_{AI}.
$$

\textbf{Case 2.} $\tau\in T(B)\setminus T(C)$.

In this case, we have $e_{k}^{(\tau)}=e_{0}^{(\tau)}$. For any $\sigma\in T(B)$, there exist $s_{\sigma}^{-1}$ and $y_{\sigma}\in\mathbb{Z}$ such that $s_{\sigma}s_{\sigma}^{-1}+m_{\sigma}y_{\sigma}=1$. We set
$$
e_{0}^{(\sigma)}\times e_{0}^{(\sigma)}=
(s_{\tau}s_{\sigma}^{-1}m_{\tau}+m_{\sigma}^2y_{\sigma})e_{0}^{(\sigma)}\, \text{if}\, \sigma\in T(B),
$$
$$
e_{i}^{(\sigma_1)}\times e_{j}^{(\sigma_2)}=0\, \text{if}\, \sigma_1\ne \sigma_2\, \text{or}\, (i,j)\ne (0,0).
$$

It follows from Remark 2.1 that these relations define a multiplication $\times$ on $G$. Then we consider the ring $(G,\times)$ as a subring of the ring $(\tilde{G},\times)$ and obtain
$$
g\times be_{0}^{(\tau)}=\dfrac{r_{\tau}}{m_{\tau}}e_{0}^{(\tau)}\times be_{0}^{(\tau)}=r_{\tau}b(s_{\tau}s_{\tau}^{-1}+m_{\tau}v_{\tau})e_{0}^{(\tau)}=r_{\tau}be_{0}^{(\tau)}\in \langle g \rangle_{AI}.
$$

Consequently,
$$
r_{\tau}A_{\tau}\subseteq \langle g \rangle_{AI}\;\text{for any }\,\tau\in T(B). \eqno (2.5)
$$ 

Now let $\tau\in T(C)$, $i\in I_{\tau}(C)$. We prove that $a_{i}^{(\tau)}be_{k}^{(\tau)}\in \langle g \rangle_{AI}$. It follows from Remark 2.1 that the relations
$$
\begin{cases}
e_{i}^{(\tau)}\times e_{i}^{(\tau)}=e_{k}^{(\tau)},\\
e_{j}^{(\sigma_1)}\times e_{s}^{(\sigma_2)}=0 
\qquad \text{if} \,(\sigma_1,\sigma_2)\ne (\tau,\tau) \, \text{or} \,(j,s)\ne (i,i),
\end{cases}
$$

define the ring $({G},\times)$. In this ring, we have
$$
g\times be_{i}^{(\tau)}=a_{i}^{(\tau)}e_{i}^{(\tau)}\times be_{i}^{(\tau)} =a_{i}^{(\tau)}be_{k}^{(\tau)}\in \langle g \rangle_{AI}.
$$

Consequently, 
$$
a_{i}^{(\tau)}A_{\tau}\subseteq \langle g \rangle_{AI}\, \text{ for all } \tau\in T(C),\, i\in I_{\tau}(C).\eqno (2.6)
$$

It follows from $(2.5)$ and $(2.6)$ that $L_{\tau}\subseteq \langle g \rangle_{AI}$ for all $\tau\in T(G)$. Therefore, 
$L=\oplus_{\tau\in T(G)} L_{\tau}\subseteq \langle g \rangle_{AI}$, whence $\langle g \rangle+L\subseteq \langle g \rangle_{AI}$.

Thus, we have proved that 
$$
\langle g \rangle_{AI}=\langle g \rangle+L.\eqno (2.7)
$$

Now we prove that if $\tau\in T(G),\, \{b_i,\,|\, i\in I_{\tau}\}\subseteq R_{\tau}$, $\ell=\text{gcd}\,(\overline{b_i}\,|\,i\in I_{\tau})$, then 
$$
\sum_{i\in I_{\tau}}b_iA_{\tau}=\ell A_{\tau}.\eqno (2.8)
$$ 

It follows from Remark 2.3 that we can represent the elements $b_i$ in the form
$b_i=\overline{b_i}c_i$, where $c_i$ is an invertible element in $R_{\tau}$, $\overline{b_i}=\ell b_i'$ and $b_i'\in\mathbb{Z}$. 

Let $x=\sum_{i\in I_{\tau}}b_ix_i\in\sum_{i\in I_{\tau}}b_iA_{\tau}$, where $x_i\in A_{\tau}$. Then
$$
x=\sum_{i\in I_{\tau}}b_ix_i=\sum_{i\in I_{\tau}}\overline{b_i}(c_ix_i)=\ell\sum_{i\in I_{\tau}}b_i'\left(c_ix_i\right)\in\ell A_{\tau}.
$$

Conversely, let $\ell x\in\ell A_{\tau}$, where $x\in A_{\tau}$. Since $\ell =\text{gcd}\,(\overline{b_i}\,|\,i\in I_{\tau})$, we can represent $\ell$ in the form $\ell =\sum_{i\in I_{\tau}}y_i\overline{b_i}$ for some $y_i\in\mathbb{Z}$. Therefore,
$$
\ell x=\left(\sum_{i\in I_{\tau}}y_i\overline{b_i}\right)x =\left(\sum_{i\in I_{\tau}}y_ib_ic_i^{-1}\right)x=
\sum_{i\in I_{\tau}}b_i(y_ic_i^{-1}x)\in \sum_{i\in I_{\tau}}b_iA_{\tau},
$$

where $c_i^{-1}$ is the element which is converse to $c_i$ in the ring $R_{\tau}$. 
Consequently, $\sum_{i\in I_{\tau}}b_iA_{\tau}=\ell A_{\tau}$.

It follows from $(2.7)$ and $(2.8)$ that $L_{\tau}=\ell_{\tau}(g)A_{\tau}$ for all $\tau\in T(G)$, whence
$$
\langle g \rangle_{AI}=\langle g\rangle +L=\langle g\rangle + \oplus_{\tau\in T(G)}L_{\tau}=\langle g\rangle + \oplus_{\tau\in T(G)}\ell_{\tau}(g)A_{\tau}.\quad \blacktriangleright
$$

\section{$afi$-Groups in the Class of $CRQ$-Groups}\label{section3}

The description of the principal absolute ideals allows us to solve the problem of describing $afi$-groups in the class $\mathcal{A}_0$. 
We recall that an \textsf{$afi$-group} is a group which has no absolute ideals other than fully invariant subgroups.

\textbf{Lemma 3.1 \cite{Pha11}.} A group $G$ is an $afi$-group if and only if $\langle g\rangle_{AI}$ is a fully invariant subgroup of the group $G$ for any $g\in G$.

\textbf{Theorem 3.2.} Any block-rigid $CRQ$-group of ring type is an $afi$-group.

$\blacktriangleleft$
Let $G$ be a block-rigid $CRQ$-group of ring type with standard representation $(2.3)$. Let $g_i\in G$,
$$
g=\sum_{\tau\in T(B)}\dfrac{r_{\tau}}{m_{\tau}}e_{0}^{(\tau)}+
\sum_{\tau\in T(C)}\sum_{i\in I_{\tau}(C)}a_{i}^{(\tau)}e_{i}^{(\tau)}\in G,\;(r_{\tau},a_{i}^{(\tau)}\in R_{\tau}).
$$
By Theorem 2.4, we have $\langle g \rangle_{AI}=\langle g \rangle +L$, where 
$$
L=\oplus_{\tau\in T(G)}L_{\tau},\quad
L_{\tau}=\ell_{\tau}(g)A_{\tau}=r_{\tau}A_{\tau}+
\sum_{i\in I_{\tau}(C)}a_{i}^{(\tau)}A_{\tau}.
$$

It is easy to see that $L$ is a fully invariant subgroup of the group $G$. Therefore, it is sufficient to prove that $\varphi(g)\in \langle g \rangle +L$ for any $\varphi\in \text{End}\,G$.

Let $\varphi\in \text{End}\,G$. Then
$$
\varphi(g)=\varphi\left(\sum_{\tau\in T(B)}\dfrac{r_{\tau}}{m_{\tau}}e_{0}^{(\tau)}+\sum_{\tau\in T(C)}\sum_{i\in I_{\tau}(C)}a_{i}^{(\tau)}e_{i}^{(\tau)}\right)=
$$
$$
=\sum_{\tau\in T(B)}\dfrac{r_{\tau}}{m_{\tau}}\varphi(e_{0}^{(\tau)})+\sum_{\tau\in T(C)}\sum_{i\in I_{\tau}(C)}a_{i}^{(\tau)}\varphi(e_{i}^{(\tau)}).\eqno(3.1)
$$
Since $\varphi(e_{i}^{(\tau)})\in A_{\tau}$ for all $\tau\in T(C)$ and $i\in I_{\tau}(C)$, we have that 
$$
\sum_{\tau\in T(C)}\sum_{i\in I_{\tau}(C)}a_{i}^{(\tau)}\varphi(e_{i}^{(\tau)})\in L.\eqno(3.2)
$$

It follows from \cite[Theorem 4.5]{Bla08} that there exists an $\alpha\in\mathbb{Z}$ such that for any $\tau\in T(B)$, the element $\varphi(e_{0}^{(\tau)})$ is of the form 
$$
\varphi(e_{0}^{(\tau)})=(\alpha+m_{\tau}y_{0}^{(\tau)})e_{0}^{(\tau)}+\sum_{i\in I_{\tau}(C)}m_{\tau}y_{i}^{(\tau)}e_{i}^{(\tau)},\;\text{where}\; y_{i}^{(\tau)}\in R_{\tau} \; (i\in I_{\tau}).
$$

Consequently,
$$
\sum_{\tau\in T(B)}\dfrac{r_{\tau}}{m_{\tau}}\varphi(e_{0}^{(\tau)})=\sum_{\tau\in T(B)}\dfrac{r_{\tau}}{m_{\tau}}(\alpha+m_{\tau}y_{0}^{(\tau)})e_{0}^{(\tau)}
+\sum_{\tau\in T(B)}\sum_{i\in I_{\tau}(C)}r_{\tau}y_{i}^{(\tau)}e_{i}^{(\tau)}=
$$
$$
=\alpha\left(\sum_{\tau\in T(B)}\dfrac{r_{\tau}}{m_{\tau}}e_{0}^{(\tau)}
+\sum_{\tau\in T(C)}\sum_{i\in I_{\tau}(C)}a_i^{(\tau)}e_{i}^{(\tau)}\right)-
\sum_{\tau\in T(C)}\sum_{i\in I_{\tau}(C)}a_i^{(\tau)}\alpha e_{i}^{(\tau)}+
$$
$$
+\sum_{\tau\in T(B)}r_{\tau}y_{0}^{(\tau)}e_{0}^{(\tau)}+
\sum_{\tau\in T(B)}\sum_{i\in I_{\tau}(C)}r_{\tau}y_i^{(\tau)}e_{i}^{(\tau)}=
$$
$$
=\alpha g-\sum_{\tau\in TCB)}\sum_{i\in I_{\tau}(C)}a_i^{(\tau)}\alpha e_{i}^{(\tau)}+\sum_{\tau\in T(B)}r_{\tau}\left(\sum_{i\in I_{\tau}}y_i^{(\tau)}e_{i}^{(\tau)}\right)\in\langle g \rangle +L.
$$ 

It follows from $(3.1)$ and $(3.2)$ that $\varphi(g)\in\langle g\rangle_{AI}$. Since $L$ is a fully invariant subgroup of the group $G$, we have that $\langle g\rangle_{AI}$ is a fully invariant subgroup of the group $G$, as well. It follows from Lemma 3.1 that $G$ is an $afi$-group.~$\blacktriangleright$

\end{document}